\documentclass{article}

\usepackage{amssymb}
\usepackage{amsmath}

\usepackage[dvips]{graphicx}

\begin{document}


\noindent \begin{center} \bf
\large{On algebraic structure of the set of prime numbers}

\end{center}

\noindent \begin{center} \bf {by: Ramin Zahedi*}

\noindent \textbf{}

\noindent \textbf{}

\noindent

\noindent 

\noindent 
\end{center}

\noindent \small

{The set of prime numbers has been analyzed, based on their algebraic and arithmetical structure. Here by obtaining a sort of linear formula for the set of prime numbers, they are redefined and identified; under a systematic procedure it has been shown that the set of prime numbers is combinations (unions and intersections) of some subsets of natural numbers, with more primary structures. In fact generally, the logical essence of obtained formula for prime numbers is similar to formula 2n - 1 for odd numbers, and so on. 
Subsequently, using obtained formula we can define all composite numbers. Finally specified examples for obtained formula are presented.}

\noindent ‫


\noindent 

\noindent 

\noindent 

\noindent \textbf{}

\noindent \textbf{1  Introduction:}

\noindent 

\noindent \textbf{}

\noindent 

\noindent

The formal definition of a prime number is as follows (using the division method):

\noindent 

\noindent 
‭
\noindent ``An integer $p{\rm >\; 1}$ is a prime if the only positive divisors of $p$ are 1 and $p$ itself.''

\noindent  

\noindent 

\noindent 

Examples of prime numbers include: 2,3,5,7,11,13,17,19,23,29,31,37...; there are infinitely prime numbers [1, 2, 37].

\noindent 

\noindent 

The prime number detection and generation has been of great interest for mathematicians all over the world for over two centuries.
Prime numbers lie at the core of some of the oldest and most perplexing questions in mathematics. Evenly divisible only by themselves and 1,
they are the building blocks of integers. In recent decades, prime numbers have emerged from their starring roles in mathematical research,
by becoming prized commodities - as elements in a cryptographic scheme widely used to keep digital messages secret [3].

\noindent 

\noindent 

Senior Max Planck Institute mathematician Don Zagier, in his article discusses prime numbers [4]: ``\dots Despite prime numbers
simple definition and role as the building blocks of the natural numbers, the prime numbers... grow like weeds among the natural numbers,
seeming to obey no other law than that of chance, and nobody can predict where the next one will sprout'' (Havil, (2003) 171, [5]).

\noindent 

Today there are many applications for primes in many scientific fields such as computer science, engineering, security science,
physics and chemistry, etc. [6-24, 27, 33-36].

\noindent 

\noindent 

In addition, there are dozens of algorithms in computer science that depend heavily on prime numbers- hashing schemes, sorting
schemes, and so on. Indirectly, as a result of studying nonlinear dynamics and chaos, Polish physicist Marek Wolf has discovered at least
two instances of fractality within the distribution of prime numbers [6-8].

\noindent 

B.L. Julia has reinterpreted the (pure mathematical) Riemann zeta function as a (thermodynamic) partition function by defining an
abstract numerical 'gas' using the prime numbers [10].

\noindent 

\noindent ------------------------------------

\noindent{\footnotesize{ *zahedi@let.hokudai.ac.jp, zahedi@logic.let.hokudai.ac.jp, zahedi.r@gmail.com, Logic Group, Hokkaido University, Japan.}}

Prime numbers are a fundamental ingredient in public-key cryptography, be it in schemes based on the hardness of factoring (e.g.
RSA), of the discrete logarithm problem, or of other computational problems. Generating appropriate prime numbers is a basic,
security-sensitive cryptographic operation [25-27, 38].

\noindent 

Complexity theory is a field in theoretical computer science, which attempts to quantify the difficulty of computational tasks and
tends to aim at generality while doing so. "Complexity" is measured by various natural computing resources, such as the amount of memory
needed, communication bandwidth, time of execution, etc. By analyzing several candidate algorithms for a problem, a most efficient one can
be easily identified; for the problem of determining the primality of an integer, the resource that could be examined, is the time of
execution [13]. There are some sorts of formulas for prime numbers; the most of these formulas have been constructed and formulated by using the
floor functions (in the field of real numbers)  [1, 2, 28-31, 37].

\noindent \textbf{}

\noindent 

\noindent The logical nature of formula for the set of prime numbers that has been obtained here (formula (21)), is similar to formula  2n - 1  
for the set of odd numbers, and so on. All these kind of algebraic formulas (also as definitions) only contain operators of the ring of integers: multiplication, 
addition and subtraction.

\noindent 

\noindent 

In fact, linear formula (21) is obtained in the same and unique process that formula 2n - 1 is formulated for the set of odd numbers. 
However above formal definition (using the division method, where the division operator is not an operator of the ring of integers) of the set of prime 
numbers seems simple, formula (21) shows clearly how this set of numbers is not a simple set.

\noindent

\noindent \textbf{} 

\noindent 

\noindent \textbf{2  A new formulation - definition:}

\noindent \textbf{}

\noindent 

\noindent 

Suppose $p_{1} ,p_{2} ,...,p_{r} $ are given prime numbers where $p_{i} $ is the \textit{ith} prime number. It follows from the
definition that in a given range of $(p_{r} ,p_{r+1}^{2} )$ any number that is not divisible by any of $p_{1} ,p_{2} ,...,p_{r} $ is a prime
number, thus let $H_{r} $ be a set of natural numbers N excluding the set of all positive multiples of $p_{1} ,p_{2} ,...,p_{r} $. 

\noindent \textbf{}

\noindent 

That is:

\noindent \textbf{}

\noindent 

\noindent \begin{center}$H_{r} =\left\{s \right. \mid s  \in {\rm N},$ and \textit s is not divisible by $\left.  p_{1} ,p_{2} ,...,p_{r}  \right\}$\[\tag{1}\]\end{center}

\noindent \textbf{} 

\noindent

\noindent 

Let $E_{1} $ be the set of natural numbers N excluding the set of all multipliers of first prime number $p_{1} =2$, define $E_{11}
$ as:

\noindent \textbf{}

\noindent 
\[E_{11} =\left\{m_{11} \right. \mid m_{11} =p_{1} x_{1} -h_{1} , \left. x_{1} \in {\rm N} \right\}          \tag{2}\]

\noindent \textbf{}

\noindent where $h_{1} =1$, we get $E_{1} =E_{11} $, and let $E_{2} $ be the set of natural numbers N excluding the set of all multipliers of
the second prime number $p_{2} =3$ , define $E_{22} $ and $E_{21} $  as:

\noindent \textbf{}

\noindent 
\[E_{21} =\left\{m_{21} \right. \mid m_{21} =p_{2} x_{2} -1, \left. x_{2} \in {\rm N} \right\},\] 
\[E_{22} =\left\{m_{22} \right. \mid m_{22} =p_{2} x_{2} -2, \left. x_{2} \in {\rm N} \right\}\]

\noindent \textbf{}

\noindent hence

\noindent \textbf{}

\noindent 
\[E_{2} =E_{21} \bigcup E_{22} =\left\{m_{2} \right. \mid m_{2} =p_{2} x_{2} -h_{2} , \left. x_{2} \in {\rm N} \right\} \tag{3}\]

\noindent \textbf{}

\noindent where $h_{2} =1,2$. Similarly let $E_{i} $be the set of natural numbers N excluding the set of all multipliers of the
\textbf{\textit{i}}\textit{th} prime number, define $E_{i1} ,E_{i2} ,...,E_{i(p_{i} -1)} $:

\noindent \textbf{}

\noindent 
\[E_{i1} =\left\{m_{i1} \right. \mid m_{i1} =p_{i} x_{i} -1, \left. x_{i} \in {\rm N} \right\},\] 
\[E_{i2} =\left\{m_{i2} \right. \mid m_{i2} =p_{i} x_{i} -2, \left. x_{i} \in {\rm N} \right\}, \]

 \begin{center}  .

  .

  .\end{center}\[E_{i(p_{i} -1)} =\left\{m_{i(p_{i} -1)} \right. \mid m_{i(p_{i} -1)} =p_{i} x_{i} -(p_{i} -1), \left. x_{i} \in {\rm N} \right\}\]

\noindent \textbf{}

\noindent then

\noindent \textbf{}

\noindent 
\[E_{i} =E_{i1} \bigcup E_{i2} \bigcup ....\bigcup E_{i(p_{i} -1)} \]

\noindent \textbf{}

\noindent 

\noindent which is equivalent to

\noindent \textbf{}

\noindent 
\[E_{i} =\left\{m_{i} \right. \mid m_{i} =p_{i} x_{i} -h_{i} , \left. x_{i} \in {\rm N} \right\}         \tag{4}\]

\noindent \textbf{}

\noindent where $h_{i} =1,2,3,...,p_{i} -1.$ It follows from the definitions above that for any set $E_{ij} $

\noindent \textbf{}

\noindent 
\[E_{ij} \bigcap E_{ik} =\emptyset        \tag{5}\]

\noindent \textbf{}

\noindent for $j\ne k$ and $i=1,2,...,r$ and $j,k=1,2,...p_{i} -1$.

\noindent \textbf{} 

\noindent 

\noindent 

It follows from (4) and (1) that

\noindent \textbf{}

\noindent 
\[H_{r} =E_{1} \bigcap E_{2} \bigcap E_{3} ....\bigcap E_{r}          \tag{6}\]

\noindent \textbf{}

\noindent 

The system of linear equations obtained from (4) and (6) define the set $H_{r} $ in natural numbers

\noindent \textbf{} 

\noindent 
\[(H)_{r} =p_{1} x_{1} -h_{1} =p_{2} x_{2} -h_{2} =p_{3} x_{3} -h_{3} =....=p_{r} x_{r} -h_{r}  \tag{7}\]

\noindent \textbf{} 

\noindent where $(H)_{r} $ is the general formula for $H_{r} $. The linear equations in (7) can be re-written as

\noindent \textbf{}

\noindent 

\noindent 
\[\left\{\begin{array}{c} {p_{1} x_{1} -h_{1} =p_{2} x_{2} -h_{2} } \\ {p_{1} x_{1} -h_{1} =p_{3} x_{3} -h_{3} } \\ {p_{1} x_{1} -h_{1}
=p_{4} x_{4} -h_{4} } \\ {.} \\ {.} \\ {.} \\ {p_{1} x_{1} -h_{1} =p_{r} x_{r} -h_{r} } \end{array}\right.\tag{8}\]

\noindent \textbf{}

\noindent for $h_{i} =1,2,3,...,p_{i} -1$ and $i=1,2,...,r$.

\noindent \textbf{}

\noindent 

\noindent 

Consider a simple linear equation in the set of  integer numbers

\noindent \textbf{}

\noindent 
\[ax-by=c      \tag{9}\]

\noindent \textbf{}

\noindent where x and y are the unknown variable and $(a,b)=1,$  $a>0,b>0$ .

\noindent 

\noindent \textbf{}

Equation (9) has infinite number of positive and negative integer solutions and in general

\noindent 
\[x=c\hat{x}'+bt,  y=c\hat{y}'+at   \tag{10}\]

\noindent where $\hat{x}'$, $\hat{y}'$are given solutions (these given solutions always exist)  of $ax'-by'=1$ and $t$ can take 
any integer value [2, 32]. Using formula (10) for the first equation in (8) we get

\noindent 
\[p_{2} x_{2} -p_{1} x_{1} =h_{2} -h_{1} =h_{2} -1\]

\noindent with the general solution of: 

\noindent \textbf{}

\noindent 
\[x_{2} =(h_{2} -1)\hat{x}'_{2} +p_{1} t_{1} ,    x_{1} =(h_{2} -1)\hat{x}'_{2} +p_{2} t_{1}          \tag{11}\]

\noindent where $\hat{x'}_{1} $, $\hat{x'}_{2} $ are given solutions for $p_{2} x'_{2} -p_{1} x'_{1} =1$ and $t_{1} $ is any integer value.

\noindent \textbf{}

Using formula (11) and the second equation of (8) we get

\noindent \textbf{}

\noindent 
\[p_{3} x_{3} -p_{1} p_{2} t_{1} =(h_{2} -1)p_{1} \hat{x}'_{1} +h_{3} -1             \tag{12}\]

\noindent \textbf{}

\noindent with the general solution of:

\noindent \textbf{}

\noindent 
\[x_{3} =[h_{3} +(h_{2} -1)p_{1} \hat{x'}_{1} -1]\hat{x'}_{3} +p_{1} p_{2} t_{2} ,\] 
\[t_{1} =[h_{3} +(h_{2} -1)p_{1} \hat{x'}_{1} -1]\hat{t}'_{1} +p_{3} t_{2}                 \tag{13}\]

\noindent where $\hat{x'}_{3} $ and $\hat{t'}_{1} $ are a given solution for  $p_{3} x'_{3} -p_{1} p_{2} t'_{1} =1$ and  $ t_{2} $ is any integer value. 

\noindent \textbf{}

Using (11), (13) and the third equation of (8) we obtain

\noindent \textbf{} 

\noindent 
\[p_{4} x_{4} -p_{1} p_{2} p_{3} t_{2} =[(h_{4} -1)+(h_{3} -1)p_{1} p_{2} \hat{t'}_{1} +(h_{2} -1)p_{1} \hat{x'}_{1} p_{3} \hat{x'}_{3} ]
\tag{14}\]

\noindent \textbf{}

\noindent 

The general solution of (14) is

\noindent \textbf{} 

\noindent 
\[x_{4} =\hat{x'}_{4} [(h_{4} -1)+(h_{3} -1)p_{1} p_{2} \hat{t'}_{1} +(h_{2} -1)p_{1} \hat{x'}_{1} p_{3} \hat{x'}_{3} ]+p_{1} p_{2} p_{3}
t_{3} ,\]\[t_{2} =\hat{t'}_{2} [(h_{4} -1)+(h_{3} -1)p_{1} p_{2} \hat{t'}_{1} +(h_{2} -1)p_{1} \hat{x'}_{1} p_{3} \hat{x'}_{3} ]+p_{4} t_{3} \tag{15}\]

\noindent \textbf{} 

\noindent 

Continuing this solution process the following general solutions are obtained:

\noindent \textbf{} 

\noindent 
\[x_{i} =\hat{x'}_{i} [(h_{i} -1)+(h_{i-1} -1)\hat{t'}_{i-3} \prod \limits _{l=1}^{i-2}p_{l} +\sum \limits _{j=2}^{i-2}[(h_{j}
-1)\hat{t'}_{j-2} \prod \limits _{l=1}^{j-1}p_{l} \prod \limits _{q=j+1}^{i-1}p_{q} \hat{x'}_{q} ]]+t_{i-1} \prod \limits _{j=1}^{i-1}p_{j} ,\]
\[ t_{i-2} =\hat{t'}_{i-2} [(h_{i} -1)+(h_{i-1} -1)\hat{t'}_{i-3} \prod \limits _{l=1}^{i-2}p_{l} +\sum \limits _{j=2}^{i-2}[(h_{j}
-1)\hat{t'}_{j-2} \prod \limits _{l=1}^{j-1}p_{l} \prod \limits _{q=j+1}^{i-1}p_{q} \hat{x'}_{q} ]]+t_{i-1} p_{i} \tag{16} \]

\noindent for $j=2,3,4,...,i-2$; $i=4,5,6,...,r$; and $\hat{t'}_{0} =\hat{x'}_{1} $. Using (16) the values of $x_{r} $ and $t_{r-2} $ can be
obtained as:

\noindent 

\noindent \textbf{}

\noindent 
\[x_{r} =\hat{x'}_{r} [(h_{r} -1)+(h_{r-1} -1)\hat{t'}_{r-3} \prod \limits _{l=1}^{r-2}p_{l} +\sum \limits _{j=2}^{r-2}[(h_{j}
-1)\hat{t'}_{j-2} \prod \limits _{l=1}^{j-1}p_{l} \prod \limits _{q=j+1}^{r-1}p_{q} \hat{x'}_{q} ]]+t_{r-1} \prod \limits _{j=1}^{r-1}p_{j}
\tag{17} \]

\noindent 
\[t_{r-2} =\hat{t'}_{r-2} [(h_{r} -1)+(h_{r-1} -1)\hat{t'}_{r-3} \prod \limits _{l=1}^{r-2}p_{l} +\sum \limits _{j=2}^{r-2}[(h_{j}
-1)\hat{t'}_{j-2} \prod \limits _{l=1}^{j-1}p_{l} \prod \limits _{q=j+1}^{r-1}p_{q} \hat{x'}_{q} ]]+t_{r-1} p_{r}\tag{18} \]

\noindent \textbf{}
\noindent \textbf{}

\noindent where $\hat{x'}_{i} $ and $\hat{t'}_{i-2} $ are given solutions of:

\noindent \textbf{}

\noindent 
\[p_{i} x'_{i} -t'_{i-2} \prod \limits _{k=1}^{i-1}p_{k} =1             \tag{19}\]

\noindent 

It is clear that given solutions  $\hat{x'}_{i} $ and $\hat{t'}_{i-2} $  always exist, as equation (19) is an especial case of equation (9).

\noindent 

Note, in (17) and (18) $t_{r-1} $ is a free integer variable. Using (17) and (18) the variable $t_{i} $ can be re-written in terms
of $t_{r-1} $. Furthermore, using (16) the variable $x_{i} $  can be re-written in terms of $t_{r-1} $, and general solutions of (19), and
$h_{i} $ and $p_{i} $. Using these, $x_{1} $ can be obtained:

\noindent 

\noindent 

\emph{\indent \[x_{1} =t_{r-1} \prod \limits _{l=2}^{r}p_{l} +(h_{r}
-1)\hat{t'}_{r-2} \prod \limits _{l=2}^{r-1}p_{l} +(h_{r-1}
-1)\hat{t'}_{r-3} p_{r}\hat{x'}_{r} \prod \limits
_{l=2}^{r-2}{p_ l}+ \sum \limits _{j=2}^{r-2}[(h_{j}
-1)\hat{t'}_{j-2} \prod \limits _{l=1}^{j-1}p_{l} \prod\limits
_{q=j+1}^{r}p_{q} \hat{x'}_{q} ] \tag{20}\]}
\indent \textbf{}

\noindent 

\noindent 

Using (20) and (7), the general formula for set $H_{r} $ (for \textit r$>$2) can be formulated as follows:

\noindent 

\noindent \textbf{} 

\noindent 

\noindent 
\[(H)_{r} =p_{1} x_{1} -1=t_{r-1} \prod \limits _{l=1}^{r}p_{l} +(h_{r} -1)\hat{t'}_{r-2} \prod \limits _{l=1}^{r-1}p_{l} +\sum \limits
_{j=2}^{r-1}[(h_{j} -1)\hat{t'}_{j-2} \prod \limits _{l=1}^{j-1}p_{l} \prod \limits _{q=j+1}^{r}p_{q} \hat{x'}_{q} ]-1 \tag{21}\]

\noindent \textbf{} 

\noindent where  $r=3,4,5,6,... $; $ j=2,3,4,...,r-1$; and  $h_{j} =1,2,3,4,...,p_{j} -1$; and parameter $t_{r-1}$ is a free integer variable, 
and $\hat{t'}_{0} =\hat{x'}_{1}$.

\noindent 

\noindent \textbf{} 

\noindent 

\textit{Following the definition of set $H_{r} $, the integer values of formula (21) are primes in the range $(p_{r} ,p_{r+1}^{2}
)$or range $[p_{r+1} ,p_{r+1}^{2} )$.}\textbf{\textit{ }}We remember that all terms in (21) are made up of prime numbers $p_{1} ,p_{2}
,p_{3} ,...,p_{r} $. Also we must note it is clear that in range \textbf{\textit{$(p_{r} ,p_{r+1}^{2} )$}}, we have at least
\textbf{\textit{$p_{r+1}^{} $}} as a prime number.\textbf{\textit{ }}Furthermore, $\hat{x'}_{i} $ and $\hat{t'}_{i-2} $  in equation (19) do not have unique values and hence formula (21) can be written in different but equivalent cases. As example, let $r=2$, using formula (11) and (7) we have:
\[(H)_{2}=p_{1} x_{1} -1 =6t_{1} +2h_{2} -3                     \tag{22} \] 

\noindent \textbf{}

\noindent prime numbers in range (3, 25) can be obtained by (22).  For $r=3$, using formula (21) (we use formula (21) for $r=3,4,5,6,7,... $)  we have

\noindent \textbf{}

\noindent 
\[(H)_{3} =30t_{2} -6h_{3} -10h_{2} +15               \tag{23}\]

\noindent \textbf{}

\noindent formula (23) gives all primes in range (5,49). For $r=4$ we may get 

\noindent \textbf{}

\noindent 

\noindent 
\[(H)_{4} =210t_{3} +90h_{4} +2184h_{3} +4550h_{2} -6825        \tag{24}\]

\noindent \textbf{}

\noindent formula (24) defines all primes in range (7,121).\textbf{ }For $r=5$ we have

\noindent \textbf{}

\noindent 

\noindent 
\[(H)_{5} =2310t_{4} -210h_{5} -18810h_{4} +114114h_{3} +190190h_{2} -285285\tag{25}\]

\noindent 

\noindent \textbf{}

\noindent formula (25) defines all primes in range (11,169). Similar formulas can be derived to obtain other primes in the proceeding
ranges.

\noindent 

\noindent 

From the definition of set $H_{r} $, it is clear that the integer values of formula (21) gives all primes and \textit{2nd} numbers, in
range $[p_{r+1}^{2} ,p_{r+1}^{3} )$; (\textit{kth} number is a number which, except itself and 1, is divisible by \textit{k} number and only
\textit{k} number of primes).

\noindent 

\noindent 

\noindent Similarly using formula (21), all primes and also all composite numbers i.e. \textit{2nd}, \textit{3rd}, \dots , \textit{kth}
numbers can be define in the range $[p_{r+1}^{k} ,p_{r+1}^{k+1} )$.

\noindent \textbf{} 

\noindent 

\noindent 

In addition using some theorems such as\textbf{ }Bertrand's postulate\textbf{ }(that states for every n $>$ 1, there is always at
least one prime \textit{p} such that n\textit{ }$<$\textit{ p }$<$ 2n), the action ranges of (21) can be expanded for larger ranges. Here
by\textit{ }Bertrand's postulate, it is easy to show that all primes and composite numbers i.e. \textit{2nd} and \textit{3rd} and \dots
\textit{k\textit th} numbers can be defined by (21), in range \textbf{$[p_{r+s}^{k} ,p_{r+s}^{k+1} )$, }if \textbf{$p_{r+1} >2^{(k+1)(s-1)} $} and
all (\textit{k}\textbf{$+1$})\textit{th} numbers i.e.:

\noindent 

\noindent 
\[p_{r+i_{1} } p_{r+i_{2} } ...p_{r+i_{k+1} } \]

\noindent are set aside from this range; where \textit{k}\textbf{$\ge $}2, \textit{s}$\ge $1, $i_{j} =$1, 2, 3, \dots , \textit{s}
and\textbf{ }$j=$1, 2, 3, \dots ,\textit{k}.

 \noindent \textbf{}

\noindent 

\noindent 

As one of the perspectives of application of formula (21), we may point to the Integer Factorisation Problem, which is a basic
discussion in cryptology and security sciences. For study of this problem, that considers the prime factors of natural numbers, we can put
equally the given number(s) to value of formula (21) and study the obtained equation(s).

\noindent 

Also formula (21), specially follow to its linear structure, can give us some new ways to study prime numbers from geometrical
points of view.

\noindent 

\noindent 

\noindent 

We believe formula (21) is a basic formula for the set of prime numbers (simply as formula 2n -1 is a  basic
formula for the set of odd numbers), and doubtless it can be useful in many fields where prime numbers are used and applied [6-24, 27, 33-36].

\noindent \textbf{}

\noindent \textbf{3  Conclusion:}

\noindent \textbf{}

As you could see in various stages of the article, in fact we assumed that prime numbers: $p_{1} ,p_{2} ,...,p_{r} $ are given,
and then through a systematic method and process the prime numbers in range of $(p_{r} ,p_{r+1}^{2} )$or range \textbf{\textit{$[p_{r+1}
,p_{r+1}^{2} )$, }}by formula (21) has been obtained. This formula is linear, and factors of this formula
only depend on $p_{1} ,p_{2} ,...,p_{r} $. The process can also be used for obtaining the next and larger ranges continually. Finally, we
could specify the set of prime numbers and define and identify them. Based on the structure of formula (21) for the set of prime numbers, we show here that the prime numbers are the result of combination of some subsets of natural numbers with more primary structure. Subsequently,
using formula (21) we define composite numbers. It should be emphasized again that the logical nature of formula (21) for the set of prime
numbers is similar to formula  2n - 1 for the set of odd numbers and so on, as it was obtained from the same (and unique) process which 
this formula is  formulated for the set odd numbers. All these kind of algebraic formulas (also as definitions) only contain operators of the ring of integers: multiplication, addition and subtraction. For more clarity, we may simply and correctly compare formula 2n -- 1 for (positive) odd numbers with formula (21) for prime numbers as follows:

\noindent \textbf{}

\noindent 

\noindent ``An odd number is a positive integer that is not divisible by 2''  (using the division method), or

\noindent ``An odd number is an integer value of algebraic linear formula 2n -- 1 in the range \textit{$[1,+\infty )$};''

\noindent \textbf{}                                                            

\noindent 

\noindent ``A prime number is a positive integer ($>$ 1) that is not divisible by any number except 1 and itself''  (using the division method), or

\noindent 

\noindent 

\noindent ``A prime number is an integer value of algebraic linear formula (21)  in the range $(p_{r} ,p_{r+1}^{2} )$ or
range \textit{$[p_{r+1} ,p_{r+1}^{2} )$}, where \textit{r} =3, 4, 5, 6,  ...,\textit{$+\infty $}''*.

\noindent \textbf{}

\noindent 

Thus formula (21) could not only be used in theoretical and logical studies of natural numbers but also especially it could be
used for practical applications of prime numbers. In addition it could  be used for study of some mathematical problems that are closed to prime
numbers, such as Riemann hypothesis.

\noindent \textbf{}

\noindent\_\_\_\_\_\_\_\_\_\_\_\_\_\_\_\_\_\_\_\_\_\_\_\_\_\_

\noindent *{\footnotesize{(there are same linear algebraic formulas for \textit{r} =1, 2;  see previous page).}

\noindent \textbf{}

\noindent \textbf{}

\noindent \textbf{References:}

\noindent \textbf{}

\noindent [1]- Crandall, R. and Pomerance, C., ``Prime Numbers,'' Springer-Verlag, New York, (2001). 

\noindent [2]- Dickson, L. E., ``History of the Theory of Numbers,'' Publisher: Amer. Mathematical  Society, Vol. 1,2,3, (1999); and Chelsea, New York, (1971).

\noindent [3]- SN: 5/25/02, p. 324, SN: 2/6/99, p. 95; Available to subscribers at:

\noindent http://www.sciencenews.org/20020525/fob4.asp.

\noindent [4]- Zaiger, D., ``The First 50 Million Prime Numbers,'' Math. Intel., (1977) 221- 224.

\noindent [5]- Havil, J., ``Gamma: Exploring Euler's Constant,'' Princeton University Press, Princeton, NJ, (2003).

\noindent [6]- M. Wolf, ``Multifractality of Prime Numbers,'' Physica A 160, (1989) 24-42.

\noindent [7]- M. Wolf, ``Random walk on the Prime Numbers,'' Physica A 250, (1998) 335-344.

\noindent [8]- M. Wolf, ``1/f noise in the distribution of Prime Numbers,'' Physica A 241, (1997) 493-499.

\noindent [9]- P. Bak, C. Tang, and K. Wiesenfeld, ``Self-organized criticality,'' Physical Review A 38, (1988) 364--374.

\noindent [10]- B.L. Julia, ``Statistical theory of numbers,'' from Number Theory and Physics (eds.  J.M. Luck, P. Moussa, and M. Waldschmidt ), Springer-Verlag, (1990).

\noindent [11]- M.C. Gutzwiller, ``Chaos in Classical and Quantum Mechanics,'' Springer-Verlag,  (1991) 307-312.

\noindent [12]- M.V. Berry and J.P. Keating, ``The Riemann zeros and eigenvalue Asymptotics,'' SIAM Review, Volume 41, No. 2, (1999) 236-266.

\noindent [13]- http://crypto.cs.mcgill.ca/\~{}stiglic/PRIMES\_P\_FAQ.html (Last updated: October 25th, 2004) and http://www.claymath.org/prizeproblems/index.htm.

\noindent [14]- http://www.maths.ex.ac.uk/\~{}mwatkins/zeta/unusual.htm.

\noindent [15]- E. Goles, O. Schulz and M. Markus, ``Prime Number Selection of Cycles In a Predator-Prey Model,'' Complexity 6 No. 4, (2001).

\noindent [16]- Joseph F. Lawler Jr., ``Identifying prime numbers with a DNA computer,'' Johns Hopkins University, Paul Ehrlich Research Award (research project), (2002). 

\noindent(See:  http://www.jhu.edu/\~{}gazette/2002/08apr02/08young.html)

\noindent [17]- J. Toh� and M.A. Soto, ``Biochemical identification of prime numbers,'' Medical Hypotheses, 53 (4), October (1999) 361-361. 

\noindent [18]- P. Dittrich, W. Banzhaf, H. Rauhe and J. Ziegler, ``Macroscopic and microscopic computation in an artificial chemistry,'' paper presented at Second German Workshops on Artificial Life (GWAL'97), Dortmund, (1997). 

\noindent [19]- R.D. Silverman, ``An Analysis of Shamir's Factoring Device'', RSA Laboratories, Bulletin, May 3 (1999). 

\noindent [20]- http://www.wolframscience.com/preview/nks\_pages/?NKS0640.gif (Stephen Wolfram explains how primes can be computed (at least theoretically) using cellular automata).

\noindent [21]- G. Mussardo, ``The quantum mechanical potential for the prime numbers,'' preprint ISAS/EP/97/153; see also R. Matthews, New Scientist, January 10th, (1998) 18.

\noindent [22]- P.W. Shor, ``Polynomial-time algorithms for prime factorization and discrete logarihtms on a quantum computer,`` SIAM J. Computing 26, (1997) 1484-1509; P.W. Shor, "Quantum computing," Documenta Mathematica Extra Volume ICM I, (1998) 467-486. 

\noindent [23]- K.Kuriyama, S.Sano and S. Furuichi, ``A precise estimation of the computational complexity in Shor's factoring algorithm,'' quant-ph/0406145.

\noindent [24]- C. Weiss, S. Page, and M. Holthaus, ``Factorising numbers with a Bose-Einstein Condensate,'' Physica A 341, (2004) 586-606.

\noindent [25]- Ueli Maurer, ``Fast Generation of Prime Numbers and Secure Public-Key Cryptographic Parameters,'' International Association for Cryptologic Research, Journal of Cryptology, Vol. 8, no. 3, (1995).

\noindent [26]- P. Beauchemin, G. Brassard, C. Crepeau, C. Goutier and C. Pomerance,``The Generation of Random Numbers that Are Probably Prime,'' International Association for Cryptologic Research, Journal of  Cryptology, Vol. 1, no. 1, (1988).

\noindent [27]- E. Kranakis, ``Primality and Cryptography,'' Series: Wiley-Teubner Series  in Computing, John Wiley and Sons Ltd, (1986).

\noindent [28]- Dudley U., ``History of Formula for Primes,'' Amer. Math. Monthly 76, (1969) 23- 28.

\noindent [29]- Guy, R. K., ``Prime Numbers'', ``Formulas for Primes'' and ``Products Taken over Primes,'' Ch. A, \S A17, and \S B48 in Unsolved Problems in Number Theory, 2${}^{nd}$ ed. New York: Springer-Verlag, (1994)  3-43, 36-41 and 102-103.

\noindent [30]- http://mathworld.wolfram.com/PrimeNumber.html.

\noindent [31]- http://mathworld.wolfram.com/PrimeFormulas.html.

\noindent [32]- Mordell, L. J., ``Diophantine Equations,'' Academic Press, New York, (1969).

\noindent [33]- Peter J. Giblin, ``Primes and Programming: Computers and Number Theory ,'' Cambridge University Press, (2004).

\noindent [34]- Hans Riesel, ``Prime Numbers and Computer Methods for Factorisation, `` (Progress in Mathematics Vol. 126), 2nd edition, Birkhauser Boston, (1994) 1- 36, 173-221, 226-237.

\noindent [35]- R. C. Vaughan, A. E. Ingham , ``The Distribution of Prime Numbers,'' Cambridge University Press, (2004).

\noindent [36]- Alan Best, ``Number Theory and Mathematical Logic: Prime Numbers  Unit 2,'' Open University Worldwide, (1996).

\noindent [37]- William Ellison, ``Prime Numbers,'' Fern Ellison, John Wiley \& Sons, (1985).

\noindent [38]- http://perltraining.com.au/\~{}jarich/Primes/.

\noindent \begin{flushleft}

\end{flushleft}

\end{document}